\documentclass[reqno]{amsart}
\usepackage{hyperref}

\AtBeginDocument{{\noindent\small
\emph{}}
\vspace{8mm}}

\begin{document}
\title[]
{On a class of nonlinear Schr$\ddot{O}$dinger equation on finite graphs}

\author[]
{Shoudong Man}

\address{Shoudong Man \newline
Department of Mathematics,
Tianjin University of Finance and Economics,
Tianjin 300222, P. R. China}
\email{manshoudong@163.com; shoudongmantj@tjufe.edu.cn}

\dedicatory{}

\thanks{The author is supported by the National Natural Science Foundation of China (Grant No. 11601368)}

\begin{abstract}
Suppose that $G=(V, E)$ is a
finite graph with the vertex set $V$ and the edge set $E$. Let $\Delta$ be the usual graph Laplacian. Consider the following nonlinear Schr$\ddot{o}$dinger type equation of the form

$$
\left \{
\begin{array}{lcr}
 -\Delta u-\alpha u=f(x,u), \\
  u\in W^{1,2}(V),\\
\end{array}
\right.
$$
on graph $G$, where $f(x,u):V\times\mathbb{R}\rightarrow\mathbb{R}$
is a nonlinear function and $\alpha$ is a parameter.
Firstly, we prove the Trudinger-Moser inequality on graph $G$, and under the assumption that $G$ satisfies the curvature-dimension type inequality $CD(m, \xi)$, we prove an integral inequality on $G$. Then by using the two inequalities, we prove that there exists a positive solution to the nonlinear Schr$\ddot{o}$dinger type equation if $\alpha<\frac{2\lambda^{2}}{m(\lambda-\xi)}$, where $\lambda$ is the eigenvalue of the graph Laplacian. Our work provides remarkable improvements to the previous results.

\end{abstract}

\maketitle
\numberwithin{equation}{section}
\newtheorem{theorem}{Theorem}[section]
\newtheorem{lemma}[theorem]{Lemma}
\newtheorem{remark}[theorem]{Remark}
\newtheorem{definition}[theorem]{Definition}
\allowdisplaybreaks

\section{Introduction}

During the past several decades,
the nonlinear Schr$\ddot{o}$dinger type equation of the form
\begin{align}\label{1.1}
 -\Delta u+b(x)u=f(x,u), u\in W^{1,2}(\Omega)
\end{align}
has been extensively studied. In equation (1.1), $\Omega\sqsubseteq R^{n},n\geq2,f(x,u):\Omega\times\mathbb{R}\rightarrow\mathbb{R}$ is a nonlinear continuous function and $b(x)\in C(\Omega,\mathbb{R})$ is a given potential. This type equation provides a good model
for developing new mathematical methods and has important applications in science and engineering. Many papers are devoted to this kind of equations
such as \cite {YYY,Y2,LYY,LN}.

 Most recently, the investigation of discrete weighted Laplacians and various equations on graphs have attracted much attention.
In \cite{AGY}, A. Grigor¡¯yan, Y. Lin and Y. Y. Yang proved that there exists a positive solution to
\begin{align}\label{1.2}
\left \{
\begin{array}{lcr}
 -\Delta u-\alpha u=|u|^{p-2}u \ \ in \ \ \Omega^{\circ}\\
  u=0\ \ on  \ \ \partial \Omega\\
\end{array}
\right.
\end{align}
on graphs for any $p>2$ if
\begin{align}\label{1.3}
\alpha<\lambda_{1}(\Omega),
\end{align}
where
\begin{align}
\lambda_{1}(\Omega)=\inf_{u \not\equiv0, u\mid\partial\Omega=0}\frac{\int_{\Omega}|\nabla u|^{2}d\mu}{\int_{\Omega}u^{2}d\mu}.
\end{align}
For more results about differential equations on graphs, we can refer to \cite{HUABIN,GLA,GAL,LIN,XA,NZL,ZD}.

In this paper, we consider a class of nonlinear Schr$\ddot{o}$dinger equation of the form
\begin{align}\label{1.5}
\left \{
\begin{array}{lcr}
 -\Delta u-\alpha u=f(x,u), \\
  u\in W^{1,2}(V)\\
\end{array}
\right.
\end{align}
on finite graph $G$. Here $f(x,u):V\times\mathbb{R}\rightarrow\mathbb{R}$
is a nonlinear function and $\alpha$ is a parameter. The equation (\ref{1.5}) can be viewed as a discrete version of (\ref{1.1}).

Firstly, we give some notations and settings.
Let $G = (V, E)$ be a graph where $V$ denotes the vertex set and $E$ denotes the edge set.
The degree of vertex $x$, denoted by $\mu(x)$, is the number of edges connected to $x$.
If $V$ contains finite vertexes, we say that $G$ is a finite graph.
 If for every vertex $x$ of $V$,  $\mu(x)$ is finite, we say that $G$ is a locally finite graph.
 When a graph is finite, this graph is certainly locally finite.
We denote $x\sim y$ if vertex $x$ is adjacent to vertex $y$. We use $(x,y)$ to denote an edge in $E$ connecting vertices $x$ and $y$.
Let $\omega_{xy}=\omega_{yx}>0$ where $\omega_{xy}$ is the edge weight.
The finite measure $\mu(x)=\sum_{y\sim x}\omega_{xy}$.
A graph $G$ is called connected if for any vertices $x,y\in V$, there exists a sequence $\{x_{i}\}_{i=0}^{n}$ that satisfies
$x=x_{0}\sim x_{1}\sim x_{2}\sim \cdot\cdot\cdot \sim x_{n}=y$.

From \cite{YLS}, for any function $u:V\rightarrow \mathbb{R}$, the $\mu$-Laplacian of $u$ is defined as
\begin{align}\label{1.6}
\Delta u(x)=\frac{1}{\mu(x)}\sum_{y\sim x}\omega_{xy}[u(y)-u(x)].
\end{align}
 The associated gradient form reads
 \begin{align}\label{1.7}
\Gamma(u,v)(x)=\frac{1}{2}\{\Delta(u(x)v(x))-u(x)\Delta v(x)-v(x)\Delta u(x)\},
\end{align}
By (\ref{1.6}) and (\ref{1.7}), we also have
 \begin{align}\label{1.9}
\Gamma(u,v)(x)=\frac{1}{2\mu(x)}\sum_{y\sim x} \omega_{xy}(u(y)-u(x))(v(y)-v(x)).
\end{align}
The length of the gradient for $u$ is
\begin{align}\label{1.10}
|\nabla u|(x)=\sqrt{2\Gamma(u,u)(x)}=(\frac{1}{\mu(x)}\sum_{y\sim x} \omega_{xy}(u(y)-u(x))^{2})^{1/2}.
\end{align}
The Ricci curvature operator on graphs $\Gamma_{2}$ by iterating $\Gamma$ as
 \begin{align}
\Gamma_{2}(u,v)(x)=&\frac{1}{2}\{\Delta\Gamma(u,v)(x)-\Gamma(u,\Delta v)(x)-\Gamma(v,\Delta u)(x)\}   \label{1.60}\\
=&\frac{1}{4}\frac{1}{d_{x}}\sum_{y\sim x}\frac{\mu_{xy}}{d_{y}}\sum_{z\sim y}\mu_{yz}[f(x)-2f(y)+f(z)]^{2}\notag\\
 &-\frac{1}{2}\frac{1}{d_{x}}\sum_{y\sim x}\mu_{xy}[f(x)-f(y)]^{2}+\frac{1}{2}[\frac{1}{d_{x}}\sum_{y\sim x}\mu_{xy}(f(x)-f(y))]^{2}.\label{1.61}
\end{align}

To compare with the Euclidean setting, we denote, for any function $u : V \rightarrow \mathbb{R}$,
\begin{align}\label{1.11}
\int_{V}ud\mu=\sum_{x\in V}\mu(x)u(x).
\end{align}

From \cite{Fan}, the eigenvalue of the Laplacian $\Delta$ on $G=(V,E)$ can be defined as
\begin{align}\label{1.12}
\lambda=\lambda_{G}=\inf \frac{\sum_{x,y\in V, y\sim x}(u(y)-u(x))^{2}\omega_{xy}}{\sum_{x\in V}u^{2}(x)\mu(x)}
=\inf\frac{\int_{\Omega}|\nabla u|^{2}d\mu}{\int_{\Omega}u^{2}d\mu}.
\end{align}
From (\ref{1.12}) we can see all eigenvalues are nonnegative and we can get $0=\lambda_{0}\leq \lambda_{1}\leq \lambda_{2}\leq\cdot\cdot\cdot\leq \lambda_{n-1}$.
For more details, we can refer to \cite{Fan}.
By (\ref{1.11}) and Lemma 1.10 in \cite{Fan},
let $u$ denote a harmonic eigenfunction achieving $\lambda$ in (\ref{1.12}),
then, for any vertex $x\in V$, we have
\begin{align}\label{1.13}
-\Delta u(x)=\frac{1}{\mu(x)}\sum_{y\sim x}\omega_{xy}[u(x)-u(y)]=\lambda u(x).
\end{align}

In \cite{YLS}, Y. Lin and S. T. Yau introduced the curvature-dimension type inequality $CD(m, \xi)$ as follows:
\begin{definition}(Curvature-dimension type inequality)
The operator $\Delta$ satisfies the curvature-dimension type inequality $CD(m, \xi)$ for some $m>1$ if for every $u$,
\begin{align}\label{1.14}
\Gamma_{2}(u,u)(x)\geq \frac{1}{m}(\Delta u(x))^{2}+\xi \Gamma(u,u)(x).
\end{align}
 We call $m$ the dimension of the operator $\Delta$ and $\xi$ the lower bound of the Ricci
curvature of the operator $\Delta$.
\end{definition}
It is easy to see that for $m < m'$, the operator $\Delta$ satisfies the curvature-dimension
type inequality $CD(m', \xi)$ if it satisfies the curvature-dimension type inequality
$CD(m, \xi)$.

\indent In \cite{YLS}, Y. Lin and S. T. Yau proved that any locally finite graph satisfies either $CD(2, \frac{2}{d}-1)$ if $d$ is finite, or $CD(2, -1)$ if $d$ is infinite, where $d=\sup_{x\in V}\sup_{y\sim x}\frac{\mu(x)}{\mu_{xy}}$.

In addition to the above curvature-dimension type inequality, we will introduce
the well-know Trudinger-Moser inequality. In \cite{MOSER}, when $p>2$, it asserts that
\begin{align}\label{1.15}
exp(\beta|u|^{\frac{p}{p-1}})\in L^{1}(\Omega).
\end{align}
Moreover, there exists a
constant $C=C(p)$ which depends only on $p $ such that
\begin{align}\label{1.16}
\sup_{||u||_{W_{0}^{1,p}(\Omega)}\leq 1}\int_{\Omega}exp(\beta
|u|^{\frac{p}{p-1}})dx\leq C(p)|\Omega|\ \ if\ \
\beta\leq\beta_{p},
\end{align}
where $p>2$ and $\beta_{p}=p
(\omega_{p-1})^{\frac{1}{p-1}},\ \ and\ \ \omega_{p-1}$ is the
measure of the unit sphere in $R^{p}$.
In the third section, we will prove the Trudinger-Moser inequality in Theorem \ref{T1} on graphs as a discrete version of (\ref{1.16}).

 Motivated by the Trudinger-Moser inequalities, we have
\begin{definition} (\cite{JM}) Suppose that $f(x,t): V\times \mathbb{R}\rightarrow \mathbb{R}$. We call that the function $f$ has
subcritical growth at $+\infty$ , if for all $\beta>0$ and $p>2$,
\begin{align}\label{1.17}
\lim_{t\rightarrow +\infty}\frac{f(x,t)}{exp(\beta|t|^{\frac{p}{p-1}})}=0.
\end{align}
\end{definition}

 In this paper, suppose that $G = (V, E)$ is a connected finite graph that satisfies the curvature-dimension type inequality $CD(m, \xi)$.
 Firstly, we prove the Trudinger-Moser inequality and an integral inequality (in Theorem \ref{T2}) on graph $G$. Then by using the two inequalities, we prove that there exists a positive solution to the nonlinear Schr$\ddot{o}$dinger type equation (\ref{1.5}) if $\alpha<\frac{2\lambda^{2}}{m(\lambda-\xi)}$,
extending the result of (\ref{1.3}) in equation (\ref{1.2}).

 Now, we state our main theorems.
\begin{theorem}\label{T1}
(Trudinger-Moser inequality on finite graphs)
Suppose that $G=(V,E)$ is a finite graph.
Then there exists a
constant $C$ which depends only on $p$ and $V$ such that
\begin{align}\label{1.18}
\sup_{||u||_{W^{1,p}(V)}\leq 1}\int_{V}exp(\beta|u|^{\frac{p}{p-1}})d\mu\leq C|V|\ \ for\ any \ \ \beta>1 \ and \ p>2,
\end{align}
where $|V|=\int_{V}d\mu(x)=$ Vol $V$,
and Vol $V$ denotes the volume of the graph $G$.
\end{theorem}

 \begin{theorem}\label{T2}
Suppose that $G = (V, E)$ is a finite graph that
satisfies the curvature-dimension type inequality $CD(m, \xi)$, and $u$ is a harmonic eigenfunction of $-\Delta$ with eigenvalue $\lambda$.
 Assume $\lambda\neq 0$. Then the following inequality holds
\begin{align}\label{1.19}
\frac{2\lambda^{2}}{m(\lambda-\xi)}\int_{V}u^{2}d\mu\leq \int_{V}|\nabla u|^{2}d\mu.
\end{align}
\end{theorem}

\begin{remark}\label{R1}
We define the norm
\begin{align}\label{R11}
||u||_{\alpha}=(\int_{V}(|\nabla u|^{2}-\alpha|u|^{2})d\mu)^{\frac{1}{2}}.
\end{align}
By using Theorem \ref{T2}, as $\alpha<\frac{2\lambda^{2}}{m(\lambda-\xi)}$, $||u||_{\alpha}$
is equivalent to the norm
\begin{align}
||u||_{W^{1,2}}=(\int_{V}|\nabla u|^{2})^{1/2}.
\end{align}

\end{remark}

 \begin{theorem}\label{T3} Suppose that $G = (V, E)$ is a finite graph that satisfies the curvature-dimension type inequality $CD(m, \xi)$.
Assume that $f: V\times \mathbb{R}\rightarrow \mathbb{R}$ satisfies the following hypotheses:\\
$(H1)$ For any $x\in V$,  $f(x,t)$ is continuous in $t\in \mathbb{R}$;\\
$(H2)$ For all $(x,t)\in V\times[0,+\infty), f(x,t)\geq 0$, and $f(x,0)=0$ for all $x\in V$;\\
$(H3)$ $f(x,t)$ has subcritical growth at $+\infty$ , i.e. $f$ satisfies (\ref{1.17});\\
$(H4)$For any $x\in V$ and $p>2$, there holds $\lim_{t\rightarrow 0+}\frac{f(x,t)}{t^{p-1}}=0$;\\
$(H5)$ There exists $q>p>2$ and $s_{0}>0$ such that if $s\geq s_{0}$, then there holds $0<qF(x,s)<f(x,s)s$ for any $x\in V$, where
$F(x,s)=\int_{0}^{s}f(x,t)dt$. \\
 Then, for any $p>2$ and
 \begin{align}\label{1.22}
 \alpha<\frac{2\lambda^{2}}{m(\lambda-\xi)},
 \end{align}
 there exists a positive solution to the equation (\ref{1.5}).
\end{theorem}

\begin{remark}\label{R2} We consider $\lambda=\lambda_{1}$ in (\ref{1.22}), where $\lambda_{1}$ is the first nonzero eigenvalue of $-\Delta$ of graph $G$.
By Lemma \ref{qulv} in the second section in this paper, we have that $\lambda_{1}\geq\frac{m\xi}{m-1}$. We can easily check that when $\frac{m\xi}{m-1}\leq\lambda_{1}<\frac{m\xi}{m-2}$,  we have $\frac{2\lambda_{1}^{2}}{m(\lambda_{1}-\xi)}>\lambda_{1}$. So we provide a remarkable improvement to the result of (\ref{1.3}) in equation (\ref{1.2}). For example, consider a connected path with two vertices a and b. It has a nonzero eigenvalue $\lambda_{1}=2$, and satisfies CD(2,1). We can check that $\frac{2\lambda_{1}^{2}}{m(\lambda_{1}-\xi)}=4>\lambda_{1}=\frac{m\xi}{m-1}=2$.

\end{remark}

 This paper is organized as follows. In Section 2, we introduce some notations and lemmas which are useful for the proof of our main theorems.
 In section 3, section 4, and section 5, we prove our main theorems.

\section{ Preliminaries }
\indent In this section, we introduce some notations and lemmas. Throughout this paper, we denote
$L^{p}(V)$ the Banach space with the norm
$||u||_{L^{p}}=(\int_{V}|u|^{p})^{1/p}$, and $W^{1,p}(V)$ can be defined as a set of all functions $u:V\rightarrow \mathbb{R}$
under the norm
\begin{align}
||u||_{W^{1,p}}=(\int_{V}|\nabla u|^{p})^{1/p}.
\end{align}

 From \cite{cyk} (See Theorem 2.1), we have the following Lemma \ref{qulv}.
\begin{lemma}\label{qulv}(\cite{cyk})
Suppose that $G = (V, E)$ is a finite graph that
satisfies the curvature-dimension type inequality $CD(m, \xi)$,
and the Ricci curvature of $G $ is at least $\xi$. Then any nonzero eigenvalue $\lambda$ of $-\Delta$ satisfies $\lambda\geq\frac{m\xi}{m-1}$.
\end{lemma}

\begin{lemma}\label{y1}
(Ambrosetti-Rabinowitz \cite{YLL}). Let $(X, ||\cdot ||)$ be a Banach space, $J \in C^{1}(X,\mathbb{R})$, $e \in X$
and $r > 0$ such that $||e|| > r$ and $b=\inf_{||u||=r}J(u)>J(0)>J(e)$.
If J satisfies the $(PS)_{c}$ condition with
 $c=\inf_{\gamma\in\Gamma}\max_{t\in[0,1]}J(\gamma(t))$,
where $\Gamma=\{\gamma\in C([0,1], X): \gamma(0)=0,\gamma(1)=e\}$,
then c is a critical value of J.
\end{lemma}

From \cite{AGY}(See Theorem 7 and Theorem 8), we have the following lemma.
\begin{lemma}\label{y2}
Let $G =(V,E)$ be a finite graph. Then for all $s>1$,  $W^{1,s}(V)$ is embedded in $ L^{q}(V)$
for all $1 \leq q \leq +\infty$. In particular, there exists a constant C depending only on $s$ and $V$
such that
\begin{align}
(\int_{V}|u|^{q}d\mu)^{1/q}\leq C(\int_{V}|\nabla u|^{s}d\mu)^{1/s}.
\end{align}
\end{lemma}
\noindent Moreover, $W^{1,s}(V)$ is pre-compact, namely, if $\{u_{k}\}$ is
bounded in $W^{1,s}(V)$, then us to a subsequence, there exists some $u\in W^{1,s}(V)$
 such that $u_{k}\rightarrow u$ in $W^{1,s}(V)$.

\section{The proof of Theorem \ref{T1}}

\begin{proof}

Let function $u$ satisfy $||u||_{W^{1,p}(V)}\leq 1$. Since $p>2$ and $\frac{p}{p-1}>1$, by Lemma \ref{y2}, we obtain
that there exist a constant $C_{0}$ such that
 \begin{align}\label{3.1}
(\int_{V}|u|^{\frac{p}{p-1}}d\mu )^{\frac{p-1}{p}}
\leq C_{0}(\int_{V}|\nabla u|^{p}d\mu )^{\frac{1}{p}}
=C_{0}||u||_{W^{1,p}(V)}
\leq C_{0}.
\end{align}
 Denote $\mu_{min}= min_{x\in V}\mu(x)$. Then (\ref{3.1}) leads to
\begin{align}\label{3.2}
||u||_{L^{\infty}(V)}\leq \frac{C_{0}}{\mu_{min}}.
\end{align}
Thus for any $\beta>1$ and $p>2$, we have
\begin{align}\label{3.3}
(\int_{V}exp(\beta|u|^{\frac{p}{p-1}})d\mu )^{\frac{p-1}{p}}
\leq exp(\frac{\beta C_{0}}{\mu_{min}})|V|^{\frac{p-1}{p}}.
\end{align}
So, we have
\begin{align}\label{3.4}
\sup_{||u||_{W^{1,p}(V)}\leq 1}\int_{V}exp(\beta|u|^{\frac{p}{p-1}})d\mu\leq C|V|,
\end{align}
where $C=(exp(\frac{\beta C_{0}}{\mu_{min}}))^{\frac{p}{p-1}}$.

\end{proof}

\section{The proof of Theorem \ref{T2}}
\begin{proof}
By (\ref{1.7}), (\ref{1.9}), (\ref{1.10}) and (\ref{1.60}), we have
\begin{align}
\Gamma_{2}(u,u)(x)=&\frac{1}{2}\{\Delta\Gamma(u,u)(x)-2\Gamma(u,\Delta u)(x)\}\notag\\
=&\frac{1}{4}\Delta|\nabla u|^{2}(x)-\Gamma(u,\Delta u)(x)\notag\\
=&\frac{1}{4}\Delta|\nabla u|^{2}(x)- \frac{1}{2}\{\Delta(u(x)\Delta u(x))-u(x)\Delta(\Delta u(x))-(\Delta u(x))^{2}\}.\label{4.1}
\end{align}
On the other hand, by (\ref{1.6}), we have
\begin{align}
&\Delta (u(x)\Delta u(x))\notag\\
&=\frac{1}{\mu(x)}\sum_{y\sim x}\omega_{xy}[u(y)\Delta u(y)-u(x)\Delta u(x)]\notag\\
&=\frac{1}{\mu(x)}\sum_{y\sim x}\omega_{xy}[u(y)\Delta u(y)-u(y)\Delta u(x)+u(y)\Delta u(x)-u(x)\Delta u(x)]\notag\\
&=\frac{1}{\mu(x)}\sum_{y\sim x}\omega_{xy}u(y)[\Delta u(y)-\Delta u(x)]+\Delta u(x)\cdot\frac{1}{\mu(x)}\sum_{y\sim x}\omega_{xy} [u(y)-u(x)]\notag\\
&=\frac{1}{\mu(x)}\sum_{y\sim x}\omega_{xy}[u(y)-u(x)][\Delta u(y)-\Delta u(x)]\notag\\
&~~ ~~~+\frac{1}{\mu(x)}\sum_{y\sim x}\omega_{xy}u(x)[\Delta u(y)-\Delta u(x)]+(\Delta u(x))^{2}\notag\\
&=\frac{1}{\mu(x)}\sum_{y\sim x}\omega_{xy}[u(y)-u(x)][\Delta u(y)-\Delta u(x)]+u(x)\Delta(\Delta u(x))+(\Delta u(x))^{2}.\label{4.2}
\end{align}
By (\ref{4.1}) and (\ref{4.2}), we have
\begin{align}\label{4.3}
\Gamma_{2}(u,u)(x)
=\frac{1}{4}\Delta|\nabla u|^{2}(x)- \frac{1}{2\mu(x)}\sum_{y\sim x}\omega_{xy}[u(y)-u(x)][\Delta u(y)-\Delta u(x)].
\end{align}
Supposing $u$ is a harmonic eigenfunction that satisfies $-\Delta u(x)=\lambda u(x)$, by (\ref{4.3}), we have
 \begin{align}
\Gamma_{2}(u,u)(x)&=\frac{1}{4}\Delta|\nabla u|^{2}(x)- \frac{1}{2\mu(x)}\sum_{y\sim x}\omega_{xy}[u(y)-u(x)][\lambda u(x)-\lambda u(y)]\notag\\
&=\frac{1}{4}\Delta|\nabla u|^{2}(x)+\frac{\lambda}{2\mu(x)}\sum_{y\sim x}\omega_{xy}[u(y)-u(x)]^{2}\notag\\
&=\frac{1}{4}\Delta|\nabla u|^{2}(x)+\frac{\lambda}{2} |\nabla u|^{2}(x).\label{4.4}
\end{align}
We consider
\begin{align}
  \sum_{x\in V}\mu(x)\Gamma_{2}(u,u)&=\frac{1}{4}\sum_{x\in V}\mu(x)\Delta|\nabla u|^{2}(x)+\frac{\lambda}{2}\sum_{x\in V}\mu(x)|\nabla u|^{2}(x)\notag\\
&=\frac{1}{4}\sum_{x\in V}\sum_{y\sim x}\omega_{xy}[ |\nabla u|^{2}(y)-|\nabla u|^{2}(x) ]+\frac{\lambda}{2}\sum_{x\in V}\mu(x)|\nabla u|^{2}(x)\notag\\
&=\frac{1}{4}\sum_{x\in V}\sum_{y\sim x}\omega_{xy}|\nabla u|^{2}(y)-\frac{1}{4}\sum_{x\in V}\sum_{y\sim x}\omega_{xy}|\nabla u|^{2}(x)
+\frac{\lambda}{2}\sum_{x\in V}\mu(x)|\nabla u|^{2}(x)\notag\\
&=\frac{1}{2}\sum_{(x,y)\in E}\omega_{xy}|\nabla u|^{2}(y)-\frac{1}{2}\sum_{(x,y)\in E}\omega_{xy}|\nabla u|^{2}(x)
+\frac{\lambda}{2}\sum_{x\in V}\mu(x)|\nabla u|^{2}(x)\notag\\
&=\frac{1}{2}\sum_{(x,y)\in E}\omega_{xy}|\nabla u|^{2}(y)-\frac{1}{2}\sum_{(y,x)\in E}\omega_{yx}|\nabla u|^{2}(y)
+\frac{\lambda}{2}\sum_{x\in V}\mu(x)|\nabla u|^{2}(x)\notag\\
&=\frac{\lambda}{2}\sum_{x\in V}\mu(x)|\nabla u|^{2}(x)\notag\\
&=\frac{\lambda}{2}\int_{V}|\nabla u|^{2}d\mu.\label{4.5}
\end{align}
On the other hand, since $G$ satisfies the curvature-dimension type inequality $CD(m, \xi)$, that is,
\begin{align}\label{4.6}
\Gamma_{2}(u,u)(x)\geq\frac{1}{m}(\Delta u(x))^{2}+ \xi\Gamma(u,u)(x),
\end{align}
 we have
\begin{align}
\sum_{x\in V}\mu(x)\Gamma_{2}(u,u)&\geq\frac{1}{m}\sum_{x\in V}\mu(x)(\Delta u(x))^{2}+ \xi\sum_{x\in V}\mu(x)\Gamma(u,u)(x)\notag\\
&=\frac{1}{m}\sum_{x\in V}\mu(x)\lambda^{2}(u(x))^{2}+ \frac{\xi}{2}\sum_{x\in V}\mu(x)\Gamma(u,u)(x)\notag\\
&=\frac{\lambda^{2}}{m}\int_{V}u^{2}d\mu+\frac{\xi}{2}\int_{V}|\nabla u|^{2}d\mu.\label{4.7}
\end{align}
By (\ref{4.5}) and (\ref{4.7}), we have
\begin{align}\label{4.8}
\frac{\lambda}{2}\int_{V}|\nabla u|^{2}d\mu\geq\frac{\lambda^{2}}{m}\int_{V}u^{2}d\mu+\frac{\xi}{2}\int_{V}|\nabla u|^{2}d\mu.
\end{align}
By (\ref{4.8}) and Lemma \ref{qulv}, since $\lambda>\xi$, we have
\begin{align}\label{4.9}
\frac{2\lambda^{2}}{m(\lambda-\xi)}\int_{V}u^{2}d\mu\leq \int_{V}|\nabla u|^{2}d\mu.
\end{align}

\end{proof}

\section{The proof of Theorem \ref{T3}}

\begin{proof}
Let $p>2$ and $\alpha<\frac{2\lambda^{2}}{m(\lambda-\xi)}$ be fixed.
 Now, we define the functional\\ $J:W^{1,2}(V)\rightarrow R$ by
\begin{align}\label{5.1}
J(u)=\frac{1}{2}\int_{V}(|\nabla u|^{2}-\alpha u^{2})d\mu-\int_{V}F(x,u^{+})d\mu,
\end{align}
where $u^{+}(x)=\max\{u(x),0\}$.  Indeed, from (H4), there exist $\tau,\delta>0$
such that if $|u|\leq \delta$ we have
\begin{align}\label{5.2}
f(x,u^{+})\leq \tau (u^{+})^{p-1}.
\end{align}
On the other hand, by (H3), there exist $c,\beta$ such that
\begin{align}\label{5.3}
f(x,u^{+})\leq
cexp(\beta|u|^{\frac{p}{p-1}}),\forall|u|\geq\delta.
\end{align}
Then we
obtain that, for $q>p$ ,
\begin{align}\label{5.4}
F(x,u^{+})\leq
cexp(\beta|u|^{\frac{p}{p-1}})|u|^{q},\forall|u|\geq\delta.
\end{align}
Combining \ref{5.2} and \ref{5.4}, we obtain that
\begin{align}\label{5.5}
F(x,u^{+})\leq \tau\frac{|u|^{p}}{p}+
c exp(\beta|u|^{\frac{p}{p-1}})|u|^{q}.
\end{align}
 By the H$\ddot{o}$lder inequality, we have
\begin{align}\label{5.6}
J(u)\geq
\frac{1}{2}||u||^{2}_{W^{1,2}}-\frac{\tau}{p}\int_{V}|u|^{p}d\mu-c(\int_{V}exp(\beta
p|u|^{\frac{p}{p-1}})d\mu)^{\frac{1}{p}}(\int_{V}|u|^{qp'}d\mu)^{\frac{1}{p'}},
\end{align}
where $\frac{1}{p}+\frac{1}{p'}=1$. By the Trudinger-Moser inequality in Theorem \ref{T1},
we obtain that
\begin{align}\label{5.7}
\int_{V}exp(\beta
p|u|^{\frac{p}{p-1}})d\mu=\int_{V}exp(\beta
p||u||_{L^{p}}^{\frac{p}{p-1}}(\frac{|u|}{||u||_{L^{p}}})^{\frac{p}{p-1}})d\mu<c|V|.
\end{align}
By Lemma \ref{y2}, there exists some constant C that depends only on p and $V$
such that
\begin{align}\label{5.8}
 (\int_{V}u^{p}d\mu)^{1/p}\leq C (\int_{V}|\nabla u|^{2}d\mu)^{1/2}.
 \end{align}
Since $q>p>2$, by (\ref{5.6}),(\ref{5.7}) and (\ref{5.8}), we can find some sufficiently small $r>0$ such that if $||u||_{W^{1,2}}=r$ we have
\begin{align}\label{5.9}
J(u)\geq
\frac{1}{2}||u||^{2}_{W^{1,2}}-C^{p}(\frac{\tau}{p}+c|V|)||u||^{p}_{W^{1,2}}.
 \end{align}
Therefore
\begin{align}\label{5.10}
\inf_{||u||_{W^{1,2}}=r} J(u)>0.
 \end{align}
By (H5), there exist two positive constants $c_{1}$ and $c_{2}$ such that
\begin{align}\label{5.11}
F(x,u^{+})\geq c_{1}(u^{+})^{q}-c_{2}.
 \end{align}
Take $u_{0}\in W^{1,2}(V)$ such that $u_{0}\geq 0$ and $u_{0}\neq0$. For any $t>0$, we have
\begin{align}\label{5.12}
J(tu_{0})\leq\frac{t^{2}}{2}||u||^{2}_{W^{1,2}}-c_{1}t^{q}\int_{V}|u_{0}|^{q}d\mu-c_{2}|V|.
\end{align}
Since $q>p>2$, we have $J(tu_{0})\rightarrow -\infty$ as $t\rightarrow +\infty$. Hence there exists some $u_{1}\in W^{1,2}(V)$
satisfying
\begin{align}\label{5.13}
J(u_{1})<0, ~~~~~ ||u_{1}||_{W^{1,2}(V)}>r.
\end{align}

Now we prove that $J(u)$ satisfies the $(PS)_{c}$ condition for any $c\in R$. To see this, we assume
$J(u_{k})\rightarrow c$ and $J'(u_{k})\rightarrow 0$ as $k\rightarrow \infty$, that is
\begin{align}\label{5.14}
\frac{1}{2}\int_{V}(|\nabla u_{k}|^{2}-\alpha u_{k}^{2})d\mu-\int_{V}F(x,u_{k}^{+})d\mu=c+o_{k}(1).
\end{align}
\begin{align}\label{5.15}
\int_{V}(|\nabla u_{k}|^{2}-\alpha u_{k}^{2})d\mu-\int_{V}u_{k}f(x,u_{k}^{+})d\mu= o_{k}(1)||u_{k}||_{W^{1,2}(V)}.
\end{align}

view of (H5), we obtain from (\ref{5.14}) and (\ref{5.15}) that $u_{k}$ is bounded in $W^{1,2}(V)$.
 Then the $(PS)_{c}$ condition follows by Lemma \ref{y2}  .
Combining (\ref{5.10}), (\ref{5.13}) and the obvious fact that $J(0)=0$, we conclude by Lemma \ref{y1} that there
exists a function $u\in W^{1,2}(V)$
such that $J(u)=\inf_{\gamma\in\Gamma}\max_{t\in[0,1]}J(\gamma(t))>0$ and $J'(u)= 0$,
where $\Gamma=\{\gamma\in C([0,1], W^{1,2}(V)): \gamma(0)=0,\gamma(1)=u_{1}\} $.
Hence there exists a nontrivial solution $u\in W^{1,2}(V)$ to the equation
\begin{align}\label{5.16}
\left \{
\begin{array}{lcr}
 -\Delta u-\alpha u=f(x,u), \\
  u\in W^{1,2}(V)\\
\end{array}
\right.
\end{align}

Testing (\ref{5.16}) by $u^{-}=\min\{u,0\}$ and noting that
\begin{align}\label{5.17}
\Gamma(u^{-},u)=\Gamma(u^{-},u^{-})+\Gamma(u^{-},u^{+})\geq |\nabla u^{-}|^{2},
\end{align}
we have
\begin{align}\label{5.18}
-\int_{V}u^{-}\Delta ud\mu-\alpha \int_{V}(u^{-})^{2}d\mu=\int_{V}u^{-}f(x,u^{+})d\mu=0.
\end{align}
This implies that $u^{-}\equiv 0$ and thus $u \geq 0$.
\end{proof}

\noindent\textbf{ Acknowledgments}
 The author thanks the referees for their comments and time.


\begin{thebibliography}{00}

\bibitem{YLL}A. Ambrosetti, P. Rabinowitz, Dual variational methods in critical point theory and applications, J. Funct. Anal. 14
(1973) 349-381.

\bibitem{cyk}Frank Bauer, Fan Chung, Yong Lin, Yuan Liu, Curvature Aspects of Graphs,
Proceedings of the American Mathematical society, Volume 145(2017) 2033-2042.

\bibitem{Fan} Fan Chung, Spectral Graph Theory, CBMS Regional Conference Series in Mathematics, American Mathematical Society, 1997.

\bibitem{cly}Fan Chung, Yong Lin, S.-T. Yau, Harnack inequalities for graphs with non-negative Ricci curvature, J. Math. Anal. Appl. 415(2014) 25-32.


\bibitem{HUABIN} HUABIN GE, A $p$-th Yamabe equation on graph, Proceedings of the American Mathematical Society, (2018)1-7.

\bibitem{AGY}
Alexander Grigoryan, Yong Lin, Yunyan Yang,Yamabe type equations on graphs,
Journal of Differential Equations, 261(9)(2016) 4924-4943.

\bibitem{GLA} Alexander Grigoryan, Yong Lin, Yunyan Yang, Kazdan-Warner equation on graph. Cal. Var. Partial Differential
Equations, 55 (4) (2016)92- 113.

\bibitem{GAL} Alexander Grigoryan, Yong Lin, Yunyan Yang, Existence of positive solutions to some nonlinear equations on
locally finite graphs. Sci. China Math., 60 (2017) 1311-1324.

\bibitem{LIN} Yong Lin, Yiting Wu, The existence and nonexistence of global solutions for a semilinear heat equation on graphs, Calculus of Variations and Partial Differential Equations, Volume 56, 2017.

\bibitem{YLS} Y. Lin, S.T. Yau, Ricci curvature and eigenvalue estimation on locally finite graphs, Math. Res. Lett. 17 (2) (2010) 345-358.

\bibitem{MOSER}J. Moser, A sharp form of an inequality by p. Trudinger,
Indiana Univ. Math. J., 20 (1970)1077-1092.

\bibitem{JM}J. M. do O, Semilinear Dirichlet problems for the p-Laplacian in R  with nonlinearities in critical
 growth range, Differential Integral Equations 9 (1996)967-979.


\bibitem{c16}Y. T. Shen, Y. X. Yao, Z. H. Chen, On a nonlinear elliptic problem with critical potential in $R^{2}$ ,
Science in China, Ser A mathematics, 47 (2004)741-755.


\bibitem{YYY} Y.Y. Yang, Existence of positive solutions to quasi-linear elliptic equations with exponential growth in the whole
Euclidean space, J. Funct. Anal., 262 (2012) 1679-1704.

\bibitem{Y2}Y.Y. Yang, L. Zhao, A class of Adams-Fontana type inequalities and related functionals on manifolds, Nonlinear
Diff. Eqn. Appl., 17 (2009) 119-135.

\bibitem{LYY} L. Zhao, Y.Y. Chang, Min-max level estimate for a singular quasilinear polyharmonic equation in $\mathbb{R}^{2m}$, J. Differential
Equations, 254 (2013) 2434-2464.

\bibitem{LN} L. Zhao, N. Zhang, Existence of solutions for a higher order Kirchhoff type problem with exponetial critical growth,
Nonlinear Anal., 132 (2016) 214-226.

\bibitem{XA} Xiaoxiao Zhang, Aijin Lin, Positive Solutions of p-th Yamabe Type Equations on Infinite Graphs,arXiv:1712.09488.

\bibitem{NZL} Ning Zhang, Liang Zhao,Convergence of ground state solutions for nonlinear Schr$\ddot{o}$dinger equations on graphs,arXiv:1705.03981.

\bibitem{ZD}Dongshuang ZHANG ,Semi-linear Elliptic Equations on Graph, J. Part. Diff. Eq.,  30(2017) 221-231.


\end{thebibliography}
\end{document}